\documentclass[12pt,leqno]{amsart}
\usepackage[latin1]{inputenc}
\usepackage[T1]{fontenc}
\usepackage [pdftex] {graphicx}
\usepackage{color}
\usepackage[francais]{babel}
\usepackage{geometry}
\usepackage{mathrsfs}
\usepackage{amsmath,amsfonts,amssymb}

\newtheorem{thm}{Théorème}[section]

\newtheorem{prop}[thm]{Proposition}
\newtheorem{defn}[thm]{Définition}
\newtheorem{rem}[thm]{Remarque}

\numberwithin{equation}{section}

\newcommand{\N}{{\mathbb{N}}}

\newcommand{\B}{{\mathbb{B}}}
\newcommand{\C}{{\mathbb{C}}}

\newcommand{\pd}{{\mathrm{d}}}
\newcommand{\h}{{\mathbb{H}}}

\newcommand{\RE}{{\mathrm{Re}}}
\newcommand{\Aut}{{\mathrm{Aut}}}

\newcommand{\dzj}{{\dfrac{\partial}{\partial{z_j}}}}
\newcommand{\dzjb}{{\dfrac{\partial}{\partial{\overline{z}_j}}}}
\newcommand{\dzn}{{\dfrac{\partial}{\partial{z_n}}}}
\newcommand{\dznb}{{\dfrac{\partial}{\partial{\overline{z}_n}}}}
\newcommand{\aut}{{\mathrm{Aut}}}
\newcommand{\un}{{\mathrm{\mathbf{1}}}}

\begin{document}
\title{Théorème de Poincaré-Alexander pour les domaines modèles.}
\author{Marianne Peyron}
\begin{abstract}
Le Théorème de Poincaré-Alexander stipule qu'une application holomorphe définie sur un ouvert de la boule unité de $\C^n$ peut, sous certaines conditions, être prolongée en un biholomorphisme de la boule unité. Dans le cadre presque complexe, la boule unité n'est plus, à biholomorphisme près, le seul domaine strictement pseudoconvexe et homogène. Un domaine strictement pseudoconvexe et homogène est biholomorphe à un ``domaine modèle``. Nous donnons dans cet article une généralisation du théorème de Poincaré-Alexander pour les domaines modèles.
\end{abstract}
\maketitle

\section*{Introduction}
Nous démontrons dans cet article le Théorème suivant :
\begin{thm} \label{PAMs}
Soit $n$ un entier strictement supérieur à $1$. Soit $(D,J)$ un domaine modèle dans $\C^n$. Soit $p\in \partial D$. Soit $U$ une boule ouverte dans $\C^n$, centrée en $p$. \\
Soit $f\, : \,  U \rightarrow \C^n$ une application de classe $\mathcal{C}^4$, pseudo-holomorphe sur $U$, telle que $f(U\cap\partial D) \subset \partial D$. Si $f$ est un difféomorphisme local en $p$, alors, $f$ s'étend en un automorphisme de $(D,J)$.
\end{thm}

Dans $\C^n$, ($n>1)$, le Théorème de Poincaré-Alexander stipule que si $U$ est une boule ouverte de $\C^n$ vérifiant $U\cap\partial\B^n\neq\varnothing$ et si $f\, : \, \overline{U\cap\B^n}\rightarrow \C^n$ est une application non constante, de classe $\mathcal{C}^1$ et holomorphe sur $U\cap\B^n$, telle que $f(\overline{U}\cap\partial \B^n) \subset \partial \B^n$, où $\B^n$ désigne la boule unité dans $\C^n$, alors, $f$ s'étend en un biholomorphisme de $\B^n$.\\
Ce Théorème a été démontré par H. Alexander dans \cite{A74} pour une application lisse, puis par S. Pin{\v{c}}uk (\cite{P75}, \cite{P77}) pour une application $\mathcal{C}^1$. (\textit{cf.} aussi \cite{R08}, th 15.3.8).\\

Dans une variété complexe, tout domaine strictement pseudoconvexe homogène est équivalent, à biholomorphisme près, à la boule unité. Ce n'est plus le cas dans une variété presque complexe : tout domaine strictement pseudoconvexe homogène est équivalent à un domaine modèle. (\textit{cf.} \cite{CGS05}, \cite{CGS06},\cite{GKK02}, \cite{L06}, \cite{L08}, \cite{BGL09}). Ainsi, il est naturel de chercher à généraliser le Théorème de Poincaré-Alexander aux domaines modèles, ce qui est fait dans le Théorème \ref{PAMs}. \\

Dans la première partie de l'article, nous définissons les notions de base telles que les structures modèles et nous rappelons un résultat concernant les jets d'ordre 2 d'une application CR, qui sera utilisé dans la deuxième partie, consacrée à la démonstration du théorème \ref{PAMs}.

\section{Préliminaires}
\subsection{Variétés presque complexes, pseudoconvexité}
Soit $M$ une variété différentielle de dimension réelle $2n$. Une\textbf{ structure presque complexe} $J$ sur $M$ est la donnée d'un isomorphisme de fibrés vectoriels différentiable, $J \,:\, T(M) \rightarrow T(M)  $ vérifiant  $J^2=-I$. On appelle \textbf{variété presque complexe} une variété différentielle $M$ munie d'une structure presque complexe $J$.\\
Une structure presque complexe $J$ est dite \textbf{intégrable} lorsque $(M,J)$ est une variété complexe. \\
Soient $(M,J)$ et $(M',J')$ deux variétés presque complexes. Une application $f\,:\,M \rightarrow M'$ de classe $\mathcal{C}^1$ est dite \textbf{$(J,J')$-holomorphe}, ou \textbf{pseudo-holomorphe}, si 
\begin{align}\label{ph}
\forall z \in M, \, \pd f_z\circ J_z=J'_{f(z)} \circ \pd f_z. 
\end{align}
Si $f\,:\,(M,J)\mapsto(M',J')$, est une bijection pseudo-holomorphe, on dira que $f$ est un $(J,J')$-biholomorphisme. On parlera de $J$-automorphisme si $(M',J')=(M,J)$. Soit $\Aut(M,J)$ le groupe des $J$-automorphismes de $M$. On dit que $(M,J)$ est \textbf{homogène} si l'action de $\Aut(M,J)$ sur $M$ est transitive.\\
Soit $\Gamma$ une sous-variété de $(M,J)$. Nous noterons  $T_{\C}\Gamma$ le complexifié de l'espace tangent de $\Gamma$: $ T_{\C}\Gamma=\C \otimes T\Gamma $. On définit l'espace tangent $J$-holomorphe de $\Gamma$ par $H^{1,0}\Gamma=\{Z \in T_{\C}\Gamma, JZ=iZ \}$.
\begin{defn}
Soient $(M,J)$ et $(N,J')$ deux variétés presque complexes. Soient $\Gamma$ et $\Gamma'$ deux sous-variétés de $M$ et $N$. Une application $\mathcal{C}^1$, $f : \Gamma\rightarrow \Gamma'$ est dite CR si $f_\star\{H^{1,0}\Gamma\} \subset H^{1,0}\Gamma'$.
\end{defn}
\begin{rem}\label{rq1}
Soit $f\,:\,(M,J)\mapsto(M',J')$ une bijection pseudo-holomorphe. Soit $\Gamma$ une sous-variété de $M$. Alors l'application $\tilde{f}=f_{|\Gamma}$ définie sur $\Gamma$ et à valeurs dans $f(\Gamma)=\Gamma'$ est CR.\\
Localement, soit $\Omega$ un ouvert de $M$ tel que $\Gamma \cap \Omega $ est non vide, soit $f\,:\,(\Omega,J)\mapsto(M',J')$ une bijection pseudo-holomorphe. Alors l'application $\tilde{f}=f_{|\Gamma\cap\Omega}$ définie sur $\Gamma\cap\Omega$ et à valeurs dans $f(\Gamma\cap\Omega)$ est CR.
\end{rem}
Soit $(M,J)$ une variété presque complexe. Si $\varphi$ est une 1-forme différentielle sur $M$, alors $J^\star\varphi$ est définie sur l'espace tangent de $M$ par $(J^\star\varphi) X=\varphi(JX)$, pour $X\in TM$. Le crochet de Lie de deux champs de vecteurs $X$ et $Y$ est le champ de vecteurs $[X,Y]$ tel que pour toute fonction $f$ de classe $\mathcal{C}^\infty$ sur $M$, on ait : $[X,Y]f=X(Yf)-Y(Xf)$.
Soit $\Gamma$ une hypersurface réelle lisse de $M$ définie par $\Gamma=\{r=0\}$. Soit $p\in\Gamma$. Le fibré tangent $J$-holomorphe de $\Gamma$ est défini par $H_p^J\Gamma=T_p\Gamma\cap JT_p\Gamma$. C'est sur cet espace qu'est définie la forme de Lévi : 
\begin{defn}\label{pseudoconvexe}
\begin{enumerate}
 \item La \textbf{forme de Lévi} de $\Gamma$ en $p$ est l'application définie sur $H_p^J\Gamma$ par $\mathcal{L}^J_\Gamma(X_p)=J^\star \pd r [X,JX]_p$, où $X$ est un champ de vecteurs de $H^J\Gamma$ tel que $X(p)=X_p$. (La définition ne dépend pas du choix d'un tel $X$).
\item Une variété presque complexe $(M,J)$ est dite \textbf{strictement $J$-pseudoconvexe} en $p$ si, pour tout $X_p$ tel que $X_p\in H_p^J\Gamma$, $\mathcal{L}^J_\Gamma(X_p)>0$.
L'hypersurface $\Gamma$ est dite \textbf{strictement $J$-pseudoconvexe} si elle l'est en tout point.
\end{enumerate}
\end{defn}


\subsection{Domaines modèles}
Nous définissons maintenant les structures modèles. Nous utilisons les notations utilisées dans \cite{CGS06}. 
\begin{defn} \label{modèle}
Une structure presque complexe $J$ sur $\C^n$ est dite \textbf{structure modèle} si $J(z)=J_{st}^{(n)}+L(z)$, où $L$ est une matrice $L=(L_{j,k})_{1 \leqslant j,k \leqslant 2n} $ telle que 
\begin{align*}
 L_{j,k}&=0 \text{ si } 1\leqslant j \leqslant 2n-2, \,1 \leqslant k \leqslant 2n,\\
 L_{j,k}&=\sum_{l=1}^{n-1}(a_l^{j,k}z_l+\overline{a}_l^{j,k}\overline{z}_l) ,\, a_l^{j,k}\in \C, \text{ si } j=2n-1,2n \text{ et } k=1,\ldots , 2n-2,
\end{align*}
où $J_{st}^{(n)}$ désigne la structure presque complexe standard sur $\C^n$.\\
La complexification $J_\C$ d'une structure modèle s'écrit comme une matrice complexe $2n\times 2n$  \\
\begin{align}\label{Jmod}
J_{\C}=\begin{pmatrix}
i&0&0&0&\ldots&0&0\\
0&-i&0&0&\ldots&0&0\\
0&0&i&0&\ldots&0&0\\
0&0&0&-i&\ldots&0&0\\
\ldots&\ldots&\ldots&\ldots&\ldots&\ldots&\ldots\\
0&\tilde{L}_{2n-1,2}(z,\overline{z})&0&\tilde{L}_{2n-1,4}(z,\overline{z})&\ldots&i&0\\
\tilde{L}_{2n,1}(z,\overline{z})&0&\tilde{L}_{2n,3}(z,\overline{z})&0&\ldots&0&-i\\
\end{pmatrix}\end{align}
avec $\tilde{L}_{2n,2i-1}(z,\overline{z})=\sum_{l=1}^{n-1}(\alpha_l^i z_l+\beta_l^i\overline{z}_l)$, où  $\alpha_l^i,\, \beta_l^i \in \C$ et $\tilde{L}_{2n,2i-1}=\overline{\tilde{L}_{2n-1,2i}}$.
\end{defn}
Une structure modèle $J$ est dite \textbf{structure modèle simple} si $\beta_l^i=0$, pour $i,l=1,\ldots,n-1$.\\
Soit $J$ une structure modèle sur $\C^n$, et $D=\{z\in \C^n,\,\RE(z_n)+P(z',\overline{z'})=0\}$, avec $P$ un polynôme homogène du second degré sur $\C^{n-1}$ à valeurs réelles. Le couple $(D,J)$ est dit \textbf{domaine modèle} si $D$ est strictement $J$-pseudoconvexe au voisinage de l'origine. Les domaines modèles sont appellés 'pseudo-Siegel domains' dans \cite{L08}.\\

Soit $\h$ le demi-plan de Siegel, défini par $\h=\{z\in\C^n,\RE(z_n)+|z'|^2=0\}$, où $z'=(z_1,\cdots,z_{n-1})$. On notera $\rho(z)=\RE(z_n)+|z'|^2$ et $\Gamma=\partial \h$.\\
Nous étudions maintenant les structures modèles et l'espace tangent $J$-holomorphe pour l'hypersurface $\Gamma$.
Soit $J_{mod}$ une structure modèle sur $\C^n$. Soit $H^{1,0}\Gamma$ l'espace tangent $J$-holomorphe . Les $(n-1)$ champs de vecteurs \begin{equation}\label{champs}
L_j=\dzj+\alpha_j(z)\dzn+\beta_j(z)\dznb,\; j=1,\ldots , n-1,
\end{equation}
 forment une base de $H^{1,0}\Gamma$, avec
\begin{equation}\label{betai}\beta_j(z) =-\dfrac{i}{2}\tilde{L}_{2n,2j-1}(z)=-\dfrac{i}{2}\sum_{l=1}^{n-1}(\alpha_l^j z_l+\beta_l^j\overline{z}_l):=\sum_{l=1}^{n-1}(a_l^j z_l+b_l^j\overline{z}_l)
 \end{equation}
et,
\begin{equation}\label{alphai}
\alpha_j(z)=2(\dfrac{i}{4}\tilde{L}_{2n,2j-1}(z)-\overline{z}_j)=-\sum_{l=1}^{n-1}(a_l^j z_l+b_l^j\overline{z}_l)-2\overline{z}_j
\end{equation}
Soit 
\begin{align}\label{T}
T=i(\dzn-\dznb).
\end{align}
 Alors, $\{T,\, L_j,\,\overline{L}_j,\, j=1,\ldots , n-1\}$ est une base de $T_{\C} \Gamma$, le complexifié de l'espace tangent de $\Gamma$.\\

Pour une matrice complexe antisymétrique $B=(b_{j,k})_{j,k=1,\ldots,n-1}$, on définit une structure presque complexe modèle simple, $J^B$, par sa complexification $J^B_{\C}$ :
\begin{align*}
 J^B_{\C}\left(\dzj\right)&=i\dzj+\sum_{k=1}^{n-1}b_{j,k}z_k\dznb,\\
 J^B_{\C}\left(\dzjb\right)&=-i\dzjb+\sum_{k=1}^{n-1}\overline{b}_{j,k}\overline{z}_k\dzn, \text{ pour } j=1,\ldots,n-1, \text{ et }\\
 J^B_{\C}\left(\dzn\right)&=i\dzn,\\
 J^B_{\C}\left(\dznb\right)&=-i\dznb.
\end{align*}

En utilisant les résultats des Théorèmes 3.4, 5.2, et 5.3 de \cite{L08}, nous obtenons :
\begin{prop}\label{modèle}
Soit $(D,J)$ un domaine modèle dans $\C^n$. Alors, il existe une matrice antisymétrique $B$ telle que $(D,J)$ est biholomorphe à $(\h,J^B)$. De plus, le biholomorphisme construit est un biholomorphisme global de $\C^n$.
\end{prop}

\subsection{Forme des applications pseudo-holomorphes locales entre deux domaines modèles}
Nous donnons ici une version locale des résultats sur la forme des applications pseudo-holomorphes entre deux domaines modèles obtenus dans \cite{L08} et \cite{BC08}. La preuve donnée dans \cite{L08} et \cite{BC08} est identique lorsque l'application est définie localement, nous ne l'incluons pas ici.
\begin{prop}\label{forme}
Soit $U$ une boule ouverte dans $\C^n$, centrée en $0$. Soient $J_s$ et $J'_s$ deux structures modèles simples et non intégrables sur $\C^n$.\\
Soit $F\, : \, U\rightarrow \C^n$ une application $(J_s,J'_s)$-holomorphe, telle que $F(U\cap\partial \h) \subset \partial \h$ et $F(0)=0$. Alors, il existe une constante réelle $c$, telle que
\begin{eqnarray}
 \forall z=(z',z_n)\in U,\, F(z,\overline{z})=(F'(z'),c z_n+\phi(\overline{z'})),
\end{eqnarray}
où $\phi$ est une application antiholomorphe (au sens standard) de $U$ dans $\C$ et $F'$ est une application holomorphe (au sens standard) de $U$ dans $\C^{n-1}$.
\end{prop}

\subsection{Jets d'ordre 2}
Nous donnons ici l'énoncé du Théorème 3.1 de \cite{P12} qui constitue un ingrédient essentiel de la preuve du Théorème \ref{PAMs}. 
\begin{thm} \label{th1}\
Soient  $J_{mod}$ et $J'_{mod}$ deux structures modèles sur $\C^n$.\\
Soit $U$ une boule ouverte centrée en $0$ dans $\C^n$. Soit $f \,:\, \partial\h\cap U \rightarrow \partial\h$ une application CR de classe $\mathcal{C}^4$. On suppose que $f$ est un difféomorphisme local en 0 et que $f(0)=0$. Alors $f$ est uniquement déterminée par ses jets d'ordre 2 en un point, et $f$ est analytique réelle.
\end{thm}
Nous mentionnons aussi le résultat obtenu par D. Zaitsev (\cite{Z12}), dans le cas où l'application CR est un difféomorphisme lisse entre deux hypersurfaces dans deux variétés presque CR de même dimension. L'hypothèse de régularité $\mathcal{C}^\infty$ de l'application ne permet pas d'utiliser ce résultat ici.
Nous signalons enfin deux résultats obtenus par F. Bertrand et L. Blanc-Centi (\cite{BBC12}), l'un en dimension réelle 4, l'autre dans le cas où la structure presque complexe est une petite déformation de la structure standard. Ici aussi, l'hypothèse sur la structure presque complexe ne permet pas d'utiliser le résultat pour la preuve du Théorème \ref{PAMs}.

\subsection{Automorphismes de $(\h,J^B)$}. \\
Les dilatations $\Lambda_\tau$ sont les automorphismes de $\h$ définis sur $\C^n$, pour $\tau>0$, par 
\begin{align}
\Lambda_\tau(z)=(\dfrac{z'}{\sqrt{\tau}},\dfrac{z_n}{\tau}).
\end{align}
Pour $\zeta\in\Gamma$, $\Psi_\zeta^B$ est l'automorphisme de $\h$ défini sur $\C^n$ par 
\begin{align}
 \Psi_\zeta^B(z)=(z'+\zeta',z_n+\zeta_n-2<z',\zeta'>_{\C}+i\mathrm{Re} B(z',\zeta')).
\end{align}
où $B(z',\zeta')=\sum_{j,k=1}^{n-1}b_{j,k}z_j\zeta_k$ et $<z',\zeta'>_{\C}=\sum_{j=1}^{n-1}z_j \overline{\zeta_j}$.\\
Nous pouvons caractériser les automorphismes de $(\h,J^B)$ à l'aide de la Proposition suivante :
\begin{prop}\label{Lee2}\cite{L08} 
Le groupe des automorphismes de $(\h,J^B)$ admet la décomposition suivante :
\begin{align*}
 \aut(\h,J^B)=\aut_{-\un}(\h,J^B)\circ \mathcal{D}\circ H^B
\end{align*}
où $\aut_{-\un}(\h,J^B)$ est le groupe d'isotropie de $-\un=(0,\ldots,0,-1)$, $\mathcal{D}=\{\Lambda_\tau,\,\tau>0 \}$ et $H^B$ est le sous-groupe de  $\aut(\h,J^B)$ défini par $H^B=\{\Psi_\zeta^B,\,\zeta\in\Gamma\}$.\\
Si $J^B$ est non intégrable, on a de plus, \[\aut_{-\un}(\h,J^B)=\{\varPhi_A(z')=(A(z'),z_n),\, A^tBA=B \text{ et } A\in U(n-1) \}.\]
\end{prop}
\begin{rem}\label{rq2}L'action de $\mathcal{D}\circ H^B$ sur $\h$ étant transitive, nous en déduisons l'homogénéité de $(\h,J^B)$. Ainsi, d'après la Proposition \ref{modèle}, tout domaine modèle $(D,J)$ est homogène.\end{rem}

\section{Démonstration du Théorème \ref{PAMs} } 
Si $n=2$, nous savons d'après \cite{CGS06} que toutes les structures modèles sont intégrables. Ainsi, dans les cas où $n=2$ ou $n$ est supérieur ou égal à $3$ et la structure $J$ est intégrable, il existe, d'après \cite{CGS06} (Prop. 2.3), un $(J,J_{st})$-biholomorphisme global de $\C^n$. Soit $\phi$ une telle application, et soit $D'=\phi(D)$. Le domaine $D'$ étant strictement pseudoconvexe et homogène (Remarque \ref{rq2}), le Théorème de Wong-Rosay (\cite{GKK02}) assure qu'il existe un biholomorphisme $\psi$ de $D$ sur $\B^n$. La version locale du Théorème de Fefferman (\cite{F76}) assure que l'on peut prolonger le biholomorphisme $\psi$ en une application lisse sur un voisinage $V$ de $p$ dans $\partial D$. Les domaines $D$ et $\B^n$ étant strictement pseudoconvexes à bords analytiques réels, la version locale du principe de réflexion de Lewy-Pinchuk (\cite{P75},\cite{L77b}), permet maintenant de prolonger l'application $\psi$ en une application holomorphe sur un voisinage $U$ de $p$ dans $\C^n$. Nous effectuons le même prolongement sur un voisinage du point $q=f(p)$, avec $q'=\phi(q)$. L'application $\psi$ se prolonge en une application holomorphe sur un voisinage $V'$ de $q'$ dans $\C^n$. \
 Quitte à rétrécir les ouverts $U$ et $V=f(U)$, nous définissons l'application $F$ sur $U''=\psi \circ \phi (U)$ par $F(z)=\psi \circ \phi \circ f \circ \phi ^{-1} \circ \psi ^{-1}(z)$. L'application $F$ vérifie les hypothèses du Théorème de Poincaré-Alexander, sur une boule contentant $p''=\psi \circ \phi (p)$ dans $U''$. Elle se prolonge donc en un biholomorphisme de $\B^n$. Par composition, nous pouvons prolonger l'application $f$ en un automorphisme de $(D,J)$.\\

Nous nous plaçons maintenant dans le cas où $n$ est supérieur ou égal à $3$ et où la structure $J$ n'est pas intégrable. Nous démontrons d'abord, grâce à la Proposition \ref{modèle}, que l'on peut se ramener au cas particulier d'un domaine modèle de la forme $(\h,J^B)$, où la structure $J^B$ est définie par une matrice antisymétrique. Nous démontrons ensuite le Théorème \ref{PAMs} dans ce cas particulier.\\

Soit $(D,J)$ un domaine modèle dans $\C^n$. Soit $p\in \partial D$. Soit $U$ une boule ouverte dans $\C^n$, centrée en $p$. 
Soit $f\, : \,  U \rightarrow \C^n$ une application de classe $\mathcal{C}^4$ et pseudo-holomorphe sur $U$, telle que $f(U\cap\partial D) \subset \partial D$, $f(p)=q$. On suppose que $f$ est un difféomorphisme local en $p$.\\
D'après la Proposition \ref{modèle}, $(D,J)$ est biholomorphe à $(\h,J^B)$. Le biholomorphisme entre $(D,J)$ et $(\h,J^B)$ est en fait un biholomorphisme global $\phi$ de $\C^n$ dans $\C^n$. Soient $p'=\phi(p)$ et $q'=\phi(q)$. Les points $p'$ et $q'$ appartiennent à $\Gamma$. De plus, $\Psi_{p'}^B(0)=p'$ et $\Psi_{q'}^B(0)=q'$. \\
On considère l'application $F\,:\,U'=\left(\Psi_{p'}^B\right)^{-1}\circ \phi(U)\mapsto \C^n$ définie par
\begin{align}
 F(z)= \left(\Psi_{q'}^B\right)^{-1}\circ\phi\circ f\circ\phi^{-1}\circ\Psi_{p'}^B(z).
\end{align}
L'application $F$ est de classe $\mathcal{C}^4$ et pseudo-holomorphe sur $U'$, vérifie $F(0)=0$ et $F(U'\cap\Gamma)\subset \Gamma$. L'application $F$ est un difféomorphisme local en $0$.\\
Nous allons montrer que l'application $F$ s'étend en un automorphisme de $(\h,J^B)$. L'application $\Psi_{q'}^B$ (resp. $\phi$) étant un $J^B$-biholomorphisme global de $\C^n$ (resp. un $(J^B,J)$-biholomorphisme global), nous obtiendrons, par composition, que l'application $f$ se prolonge en un automorphisme de $(D,J)$.\\

Nous allons d'abord calculer le jet d'ordre $2$ en $0$ de l'application $F$. Nous verrons ensuite que ce jet d'ordre $2$ en $0$ coïncide avec le jet d'ordre $2$ en $0$ d'un automorphisme de $(\h,J^B)$. Nous obtiendrons le résultat recherché en utilisant le Théorème \ref{th1}.\\
La structure $J$ n'étant pas intégrable, d'après la Proposition \ref{forme}, il existe une constante réelle $c$, telle que
\begin{align}\label{forme'}
 \forall z=(z',z_n)\in U,\, F(z,\overline{z})=(F'(z'),c z_n+\phi(\overline{z'})),
\end{align}
où $\phi$ est une application antiholomorphe (au sens standard) de $U$ dans $\C$ et $F'$ est une application holomorphe (au sens standard) de $U$ dans $\C^{n-1}$. Puisque $F(0)=0$, le développement limité de $F$ en $0$ est de la forme 
\begin{align}
F_j(z)&= a^j_1 z_1+\ldots +a^j_{n-1}z_{n-1}+\sum_{k,l=1}^{n-1}a^j_{k,l}z_kz_l+o(|z|^2)\text{ pour }j=1,\ldots,n-1,\label{F_j}\\
F_n(z,\overline{z})&=c z_n+\sum_{k=1}^{n-1}a^n_{\overline{k}}\overline{z}_k+\sum_{k,l=1}^{n-1}a^n_{\overline{k},\overline{l}}\overline{z}_k\overline{z}_l+o(|z|^2).\label{F_n}
\end{align}

Puisque $F$ vérifie les équations de pseudo-holomorphie (\ref{ph}), on obtient, d'après (\ref{Jmod}), pour $j=1,\ldots,n-1$, pour $z\in U'$,
\begin{align}\label{holo}
 \sum_{l=1}^{n-1}\tilde{L}_{2n,2l-1}(F(z))\dfrac{\partial F_l}{\partial z_j}(z)=2i\dfrac{\partial \overline{F_n}}{\partial z_j}(z)+\tilde{L}_{2n,2j-1}(z)\dfrac{\partial \overline{F_n}}{\partial \overline{z_n}}(z),
\end{align}
avec $\tilde{L}_{2n,2l-1}(w)=\sum_{k=1}^{n-1}b_{l,k}w_k $. D'où, d'après (\ref{F_n}),
\begin{align}\label{holobis}
 \sum_{l=1}^{n-1}\left(\sum_{k=1}^{n-1}b_{l,k}F_k(z)\right)\dfrac{\partial F_l}{\partial z_j}(z)=2i\dfrac{\partial \overline{F_n}}{\partial z_j}(z)+c\left(\sum_{k=1}^{n-1}b_{j,k}z_k\right).
\end{align}
Puisque $F(0)=0$, les termes constants sont nuls dans le membre de gauche de (\ref{holobis}). On obtient donc $a^n_{\overline{j}}=0$, pour $j=1,\ldots,n-1$. \\
De plus, puisque $F(U'\cap\Gamma) \subset \Gamma$, on a, pour $z\in\Gamma\cap U'$, 
\begin{align}\label{gamma}
\rho(F(z))=0.
 \end{align}
C'est à dire, pour $z\in\Gamma\cap U'$, 
\begin{align}\label{gammabis}
\RE(F_n(z,\overline{z}))+\sum_{j=1}^{n-1}F_j(z)\overline{F_j(z)}=0.
 \end{align}
D'après (\ref{F_j}) et (\ref{F_n}), pour $k,l=1,\ldots n-1$, les seuls termes en $\overline{z}_k\overline{z}_l$ dans le membre de gauche de l'égalité (\ref{gammabis}) sont $a^n_{\overline{k},\overline{l}}\overline{z}_k\overline{z}_l$. On obtient donc $a^n_{\overline{k},\overline{l}}=0$.\\
Ecrivons le développement limité de $F_n$ à l'ordre 3 : 
\begin{align*}
F_n(z,\overline{z})&=c z_n+\sum_{p,k,l=1}^{n-1}a^n_{\overline{p},\overline{k},\overline{l}}\overline{z}_p\overline{z}_k\overline{z}_l+o(|z|^3).
\end{align*}
L'égalité (\ref{gammabis}) donne, pour $z\in\Gamma \cap U'$ :
\begin{align}\label{gammater}
&c\RE(z_n)+\RE\left(\sum_{p,k,l=1}^{n-1}a^n_{\overline{p},\overline{k},\overline{l}}\overline{z}_p\overline{z}_k\overline{z}_l+o(|z|^3)\right)\\
&+\sum_{j=1}^{n-1}\left(\sum_{l=1}^{n-1}a^j_l z_l+\sum_{k,l=1}^{n-1}a^j_{k,l}z_kz_l+o(|z|^2)\right)\left(\sum_{l=1}^{n-1}\overline{a^j_l} \overline{z_l}+\sum_{k,l=1}^{n-1}\overline{a^j_{k,l}}\overline{z_k}\overline{z_l}+o(|z|^2)\right)=0.\notag
\end{align}
En écrivant le terme en $\overline{z}_kz_pz_l$ dans l'équation (\ref{gammater}), pour $k,p,l=1,\ldots,n-1$, on obtient :
\begin{align}\label{s2}
 a^1_{p,l}\overline{a^1_k}+a^2_{p,l}\overline{a^2_k}+\ldots+a^{n-1}_{p,l}\overline{a^{n-1}_k}=0.
\end{align}
Etant donné la forme de $F$ donnée en (\ref{F_j}) et (\ref{F_n}), $F$ étant un difféomorphisme local en $0$, la matrice $A=(a^i_k)_{i,k=1,\ldots,n-1}$ est inversible. En inversant le système composé des équations (\ref{s2}) pour $k=1,\ldots,n-1$, à $p$ et $l$ fixés, nous obtenons donc $a^j_{p,l}=0$, pour $j,p,l=1,\ldots,n-1$.
Le développement limité de $F$ en $0$ est donc de la forme 
\begin{align*}
 F_j(z,\overline{z})&= a^j_1 z_1+\ldots +a^j_{n-1}z_{n-1}+o(|z|^2)\text{ pour }j=1,\ldots,n-1,\\
F_n(z,\overline{z})&=c z_n+o(|z|^2).
\end{align*}
L'égalité (\ref{gammater}) devient, pour $z\in\Gamma\cap U'$ :
\begin{align}\label{gammaqua}
&c\RE(z_n)+\RE\left(\sum_{p,k,l=1}^{n-1}a^n_{\overline{p},\overline{k},\overline{l}}\overline{z}_p\overline{z}_k\overline{z}_l+o(|z|^3)\right)\\
&+\sum_{j=1}^{n-1}\left(\sum_{l=1}^{n-1}a^j_l z_l+o(|z|^2)\right)\left(\sum_{l=1}^{n-1}\overline{a^j_l} \overline{z_l}+o(|z|^2)\right)=0.\notag
\end{align}
De plus, si $z\in \Gamma\cap U'$, $\RE(z_n)=-|z'|^2$. Lorsque nous écrivons l'égalité des coefficients des termes en $z_p \overline{z}_p$ dans l'égalité (\ref{gammaqua}), nous obtenons, pour $p=1,\ldots,n-1$,
\begin{align}\label{pp}
 -c+\sum_{j=1}^{n-1}a^j_p\overline{a^j_p}=0.
\end{align}
Lorsque nous écrivons l'égalité des coefficients des termes en $z_p\overline{z}_k$ dans l'égalité (\ref{gammaqua}), nous obtenons, pour $p,k=1,\ldots,n-1$, $p\neq k$,
\begin{align}\label{pk}
\sum_{j=1}^{n-1}a^j_p\overline{a^j_k}=0.
\end{align}
Les équations (\ref{pp}) et (\ref{pk}) signifient que $A^t\overline{A}=cI_{n-1}$.\\
Définissons l'application $G=(G_1,\ldots,G_n)$ sur $\C^n$ par 
\begin{align}
 G_j(z,\overline{z})&= a^j_1 z_1+\ldots +a^j_{n-1}z_{n-1}\text{ pour }j=1,\ldots,n-1,\label{G1}\\
G_n(z,\overline{z})&=c z_n.\label{G2}
\end{align}
En posant $A'=({a'}_{i,j})_{i,j=1,\ldots,n-1}=\left(\dfrac{a_{i,j}}{\sqrt{c}}\right)_{i,j=1,\ldots,n-1}$, on a $G(z)=\varPhi_{A'}\circ\Lambda_{\frac{1}{c}}\circ\Psi_0^B$. 
Ainsi, d'après la Proposition \ref{Lee2}, l'application $G$ est un $J^B$-automorphisme de $\h$. Son jet d'ordre $2$ en $0$ coïncide avec celui de $F$. D'après la Remarque \ref{rq1}, l'application $\tilde{F}$, restriction de l'application $F$ à $\Gamma \cap U'$ est CR, vérifie les hypothèses du Théorème \ref{th1} : l'application $\tilde{F}$ est donc uniquement déterminée par son jet d'ordre $2$ en $0$. Les applications $\tilde{F}$ et $G$ sont donc égales sur $U'\cap\Gamma$. Ainsi, l'application $\tilde{F}$ se prolonge en un automorphisme de $(\h,J^B)$.\\

Montrons maintenant que l'application $F$ se prolonge en un automorphisme de $(\h,J^B)$. D'après l'écriture (\ref{forme'}), l'application $F$ est analytique réelle sur $U'$. Nous allons démontrer que ses dérivées en $0$ coïncident avec celles de l'application $G$, qui est aussi analytique réelle. Les applications $F$ et $G$ seront donc égales sur $U'$, et nous obtiendrons que l'application $F$ se prolonge en un automorphisme de $(\h,J^B)$.\\
D'après l'écriture (\ref{forme'}) et la définition de l'application $G$ donnée dans (\ref{G1}) et (\ref{G2}), nous obtenons directement que :
\begin{align*}
 \dfrac{\partial^k F^j}{\partial z^\alpha \partial \overline{z}^\beta}(0)=\dfrac{\partial^k G^j}{\partial z^\alpha \partial \overline{z}^\beta}(0)&\text{ pour }j=1,\ldots,n-1,\; k\geq 2,\, \alpha=(\alpha_1,\ldots,\alpha_n),\,\beta=(\beta_1,\ldots,\beta_n),\,\\
&|\alpha|+|\beta|=k,\,\text{ et } \alpha_n\neq0 \text{ ou l'un des } \beta_i\neq 0.\\
 \dfrac{\partial^k F^n}{\partial z^\alpha \partial \overline{z}^\beta}(0)=\dfrac{\partial^k G^n}{\partial z^\alpha \partial \overline{z}^\beta}(0)&\text{ pour }j=1,\ldots,n-1,\; k\geq 2,\, \alpha=(\alpha_1,\ldots,\alpha_n),\,\beta=(\beta_1,\ldots,\beta_n),\,\\
&|\alpha|+|\beta|=k,\,\text{ et } \beta_n\neq0 \text{ ou l'un des } \alpha_i\neq 0.\\
\end{align*}
Nous rappelons que les champs $\{T,\, L_k,\,\overline{L}_k,\, k=1,\ldots , n-1\}$ définis en (\ref{champs}) et (\ref{T}) constituent une base du complexifié de l'espace tangent $T_{\C} \Gamma$. Les applications $F$ et $G$ étant égales sur $U'\cap\Gamma$, nous avons, pour tout $p \in U'\cap\Gamma $, pour tout $j=1,\ldots,n$, pour tout $t\in\N$, pour tous multi-indices $\alpha=(\alpha_1,\ldots,\alpha_n)$ et $\beta=(\beta_1,\ldots,\beta_n)$, 
\begin{align*}
T^t_p L^\alpha_p \overline{L}^\beta_p F_j(p)=0.
 \end{align*}
En particulier, pour $p=0$, $j=1,\ldots,n-1$, $t=0$, $\alpha=(\alpha_1,\ldots,\alpha_{n-1},0)$ et $\beta=(0,\ldots,0)$, on obtient,
\begin{align*}
 \dfrac{\partial^k F^j}{\partial z^\alpha }(0)=\dfrac{\partial^k G^j}{\partial z^\alpha }(0)&\text{ pour }j=1,\ldots,n-1,\; k\geq 2,\, \alpha=(\alpha_1,\ldots,\alpha_{n-1},0),\,|\alpha|=k.
\end{align*}
De même, pour $p=0$, $j=n$, $t=0$, $\alpha=(0,\ldots,0)$ et $\beta=(\beta_1,\ldots,\beta_{n-1},0)$, on obtient,
\begin{align*}
 \dfrac{\partial^k F^n}{ \partial \overline{z}^\beta}(0)=\dfrac{\partial^k G^n}{ \partial \overline{z}^\beta}(0)&\text{ pour }j=1,\ldots,n-1,\; k\geq 2,\, \,\beta=(\beta_1,\ldots,\beta_{n-1},0),\,|\beta|=k.
\end{align*}
Ainsi, toutes les dérivées partielles en $0$ de l'application $F$ coïncident avec celles de l'application $G$, qui est aussi analytique réelle. Les applications $F$ et $G$ sont donc égales sur $U'$, et l'application $F$ se prolonge en un automorphisme de $(\h,J^B)$. Ceci termine la démonstration.\\

\providecommand{\bysame}{\leavevmode\hbox to3em{\hrulefill}\thinspace}
\providecommand{\MR}{\relax\ifhmode\unskip\space\fi MR }
\providecommand{\MRhref}[2]{%
  \href{http://www.ams.org/mathscinet-getitem?mr=#1}{#2}
}
\providecommand{\href}[2]{#2}

\vskip 0,5cm
{\small
\noindent Marianne Peyron\\
(1) UJF-Grenoble 1, Institut Fourier, Grenoble, F-38402, France\\
(2) CNRS UMR5582, Institut Fourier, Grenoble, F-38041, France\\
{\sl E-mail address} : marianne.peyron@ujf-grenoble.fr
}
\end{document}